\documentclass[12pt]{article}
\usepackage{amsmath,amsthm,amssymb}
\newtheorem{theorem}{Theorem}

\newcommand{\xtw}{{\tilde x}}

\newcommand{\beq}{\begin{equation}}
\newcommand{\eeq}{\end{equation}}

\begin{document}

\title{ \large The Lagrange Inversion Theorem with Limited Smoothness }

\author{
{\normalsize J. Thomas Beale} \\
{\normalsize{\em Department of Mathematics, Duke University, Box 90320}}\\
{\normalsize{\em Durham, North Carolina 27708, U.S.A.}}\\
}

\date{}

\maketitle

\begin{abstract}  The Lagrange Inversion Theorem provides a formula for the power series of the inverse of an analytic function.  We present a straightforward extension to functions with a finite number of derivatives.  
\end{abstract}

The Lagrange Inversion Theorem provides a formula for the power series of the inverse of an analytic function.  The book [3] has a historical discussion and a proof using the Cauchy integral formula. The Theorem is stated in various references, e.g. [1,5].
 It is useful to be able to apply the formula more generally to functions of limited smoothness.
The purpose of this note is to provide a straightforward extension to functions with a finite number of derivatives, since we have not been able to find such a version available.  Related results extending a different version to nonanalytic functions are given in [2,4], and there are results for formal power series.

The Lagrange Theorem is stated and proved in [3] as Theorem 2.3.1.  We begin with a paraphrase of that statement:

\begin{theorem}
Suppose $f$ is an analytic function of $z \in \mathbb{C}$ in a neighborhood of 
$z = a$,
$f(a) = b$, and $f'(a) \neq 0$.  Then there are sufficiently small open disks $B_1$ about $a$ and $B_2$ about $b$ so that for each $w \in B_2$ there is a unique root $z \in B_1$ of $f(z) = w$.  Set $g(w) = z$.  Then the inverse function $g$ is analytic and 
\beq g(w)  = a + \sum_{n=1}^\infty \frac{(w-b)^n}{n!} D^{n-1} \left[\phi(z)^n \right]_{z=a} \eeq
where $D = d/dz$ and $\phi$ is the analytic function
\beq \phi(z) = \frac{z-a}{f(z)-b}\,,\quad z \neq a \,, \qquad
\phi(a) = \frac{1}{f'(a)}  \eeq
\end{theorem}

To obtain this statement from Thm. 2.3.1 in [3], in (2.20) we define $\phi$ as above, choose $t = w - b$, and set $\psi(z) = z$.  Then for $z \neq a$ with $z$ near $a$, (2.20) in [3] is equivalent to $w = f(z)$, whereas for $z = a$ it is $w = b = f(a)$.
The formula (1)-(2) is equivalent to those in [1] and the NIST Digital Library of Mathematical Functions.
We will use this theorem, as well as the usual Inverse Function Theorem, to verify the following version for real-valued functions that are smooth rather than analytic:

\begin{theorem}
Suppose that $f$ is a real-valued $C^{N+1}$ function of $x \in \mathbb{R}$, 
with $N \geq 1$, defined near
$x = a$, with $f(a) = b$ and $f'(a) \neq 0$.  Then there are open intervals $I$ about $a$ and $J$ about $b$, sufficiently small, so that for each $y \in J$ there is a unique $x \in I$ with $f(x) = y$.  Set $g(y) = x$.  Then $f: g(J) \to J$ and $g: J \to g(J)$ are inverses and
\beq g(y)  = a + \sum_{n=1}^N \frac{(y-b)^n}{n!} D^{n-1} \left[\phi(x)^n \right]_{x=a} \,+\, O(|y-b|^{N+1})          \eeq
where $D = d/dx$ and $\phi$ is the $C^N$ function
\beq \phi(x) = \frac{x-a}{f(x)-b}\,, \quad x \neq a\,, 
           \qquad \phi(a) = \frac{1}{f'(a)}  \eeq
\end{theorem}

\smallskip

To begin the proof, we may assume that $a = 0$ and $b = 0$, since we can translate $x$ by $a$ and $f$ by $b$.  With $f \in C^{N+1}$ and $f(0) = 0$, it is evident
from $f(x)/x = \int_0^1 f'(tx)\,dt $
that $f(x)/x$ is $C^N$, and since $f'(0) \neq 0$,
the same is true for $\phi(x) = x/f(x)$.
Now let $p$ be the Taylor polynomial
of order $N$, so that $f(x) = p(x) + O(|x|^{N+1})$.  We will apply the above theorem to $p$, since it is analytic.  By the usual Inverse Function Theorem, $f$ and $p$ both have $C^{N+1}$ inverses for $y$ in a sufficiently small interval $J$. 
From Theorem 1 we have
\beq p^{-1}(y) = \sum_{n=1}^N \frac{y^n}{n!} D^{n-1} \left[\psi(x)^n \right]_{x=0} + O(|y|^{N+1})\,, \qquad \psi(x) = \frac{x}{p(x)} \eeq
and our task is to verify that
\beq f^{-1}(y) = \sum_{n=1}^N \frac{y^n}{n!} D^{n-1} \left[\phi(x)^n \right]_{x=0} + O(|y|^{N+1})\,, \qquad \phi(x) = \frac{x}{f(x)} \eeq

We first show that
\beq f^{-1}(y) - p^{-1}(y) = O(|y|^{N+1}) \quad \mbox{as} \quad y \to 0 \eeq
 Given $y \in J$,
let $x = f^{-1}(y)$ and $\xtw = p^{-1}(y)$.  Then
\beq 0 = f(x) - p(\xtw) = f(x) - f(\xtw) + O(|\xtw|^{N+1})\eeq
Since $x = 0$ at $y = 0$, $|x| \leq M|y|$ and $|\xtw| \leq M|y|$ for $y$ in some interval about $0$ and some $M$.  Then for $|y|$ small enough, $|x|$ and $|\xtw|$ are small, so that $|f'(\xi)| \geq |f'(0)|/2$ for $\xi$ between $x$ and $\xtw$.  Then
$|f(x) - f(\xtw)| \geq |f'(0)| |x - \xtw|/2$.  Since $f'(0) \neq 0$,
it follows from (8) that
$|x - \xtw| = O(|\xtw|^{N+1})$ and also $|x - \xtw| = O(|y|^{N+1})$
since $\xtw$ is bounded by $y$. Thus (7) holds.
% As a consequence the Taylor polynomials for $f^{-1}$ and $p^{-1}$ have the same % coefficients for $n \leq N$.

The functions $f(x)/x$ and $p(x)/x$ are $C^N$ and have the same Taylor expansion to $N-1$, with remainder $O(|x|^N)$.  The same is true for $\phi(x) = x/f(x)$ and $\psi(x) = x/p(x)$, according to standard formulas for series of reciprocals.
% in the DLMF, 1.9.53 to1.9.56.  
When we multiply series, the coefficient of $x^n$ in the product depends only on the coefficients of $x^m$ for $m \leq n$ in the factors.  Thus, since the expansions for  $\phi(x)$ and $\psi(x)$ match to order $N-1$, the same is true for 
$\phi(x)^n$ and $\psi(x)^n$.  Then $\phi(x)^n - \psi(x)^n$ is $O(|x|^N)$, and its Taylor polynomial of order $N-1$ is zero.  Therefore $D^{n-1}\phi^n -
D^{n-1}\psi^n =  0$ at $x=0$ for $n \leq N$.   Together with (7) we have now  verified that (6) is equivalent to (5).

\bigskip

{\bf References}

\smallskip

1.  G. E. Andrews, R. Askey, and R. Roy, {\em Special Functions}, Encyclopedia of Mathematics and its Applications, Vol. 71, Cambridge Univ. Press, 1999.

\smallskip

2. N. Grossman, {\em A $C^\infty$ Lagrange inversion theorem}, Amer. Math. Monthly 112 (2005), 512–514.

\smallskip

3. S. G. Krantz and H. R. Parks, {\em The Implicit Function Theorem, History, Theory, and Applications},
Birkhauser, 2002.

\smallskip

4. S. G. Krantz and H. R. Parks, {\em The Lagrange inversion theorem in the smooth case},
J. Math. Anal. Appl. 340 (2008), 1263–1270.

\smallskip

5. E. T. Whittaker and G. N. Watson, {\em A Course of Modern Analysis}, 4th ed., Cambridge Univ. Press, 1927.

\end{document}